\documentclass[11pt]{article}
\usepackage{amsfonts,amssymb}
\usepackage{graphicx}
\usepackage[margin=35mm]{geometry}

\newtheorem{LEM}{Lemma}
\newtheorem{PROP}{Proposition}
\newtheorem{OBS}{Observation}
\newtheorem{COR}{Corollary}

\newtheorem{CONJ}{Conjecture}
\newenvironment{Proof}{\textbf{Proof.} }{$\Box$\smallskip}

\newcommand{\TD}{\mathop{\mathit{TD}}}
\newcommand{\DEF}{\mathop{\mathit{df}}}
\newcommand{\MULTI}[2]{K^{(#1)}_#2}

\newcommand{\dd}[1]{\textbf{#1}}

\begin{document}

\begin{center}
\textbf{\Large Colorings with Fractional Defect} \smallskip \\
Wayne Goddard and Honghai Xu \\
Department of Mathematical Sciences, Clemson University
\end{center}

\footnotetext{Corresponding author:
Wayne Goddard, goddard@clemson.edu, phone +1-864-656-0186, fax +1-864-656-0145}

\begin{abstract}
Consider a
coloring of a graph such that each vertex is assigned a fraction of
each color, with the total amount of colors at each vertex summing to~$1$.
We define the fractional defect of a vertex $v$ to be the sum
of the overlaps with each neighbor of~$v$, and
the fractional defect of the graph to be the maximum of the defects over all vertices.
Note that this coincides with the usual definition of defect  if every vertex is monochromatic.
We provide results on the minimum fractional defect of $2$-colorings of some graphs.
\end{abstract}

\section{Introduction}

In a usual vertex coloring of a graph, every vertex  is
assigned one color and that color is different from each of its
neighbors. We consider here a two-fold generalization of this:
a vertex can receive multiple colors and can overlap slightly with each neighbor.

Specifically, each vertex is assigned a fraction of
each color, with the total amount of colors at each vertex summing to~$1$.
The \dd{(fractional) defect} of a vertex $v$ is defined to be the sum
of the overlaps over all colors and all neighbors of $v$.
For example, if vertex~$v$ receives $\frac{1}{3}$ red and $\frac{2}{3}$ blue,
while its neighbor $w$ receives $\frac{2}{5}$ red, $\frac{2}{5}$ blue, and $\frac{1}{5}$ white,
then $w$ contributes $\frac{1}{3}+\frac{2}{5}$ to the defect of $v$, and vice versa.
We define
the \dd{(fractional) defect} of the graph as the maximum of the defects over all vertices.
Note that if every vertex is monochromatic (has only one color), then our fractional defect
coincides with the usual definition of defect (see for example~\cite{CCW86defective});
and that defective colorings are also called improper colorings.

The idea of assigning vertices multiple colors has been used most
notably in fractional colorings (e.g.~\cite{PU02planar,LZ04fracDistance}), but also for example in
$t$-tone colorings~\cite{CKK13tTone}.
Like in $t$-tone colorings (and unlike in fractional colorings), we
consider here the situation where one pays for each color used, regardless
of how much the color is used.
Note that for proper colorings, allowing one to color a vertex with multiple colors
does not yield anything new. For, one can just choose for each
vertex $v$ one color present at $v$ and recolor it entirely that color,
and therefore the minimum number of colors needed is just the chromatic number.
Similarly, with the usual definition of the defect of a vertex as the number of
neighbors that share a color, there is no advantage to using more than one color
at a vertex. But we consider colorings where a vertex overlaps
only slightly with each neighbor.

\begin{figure}[!h]
\centerline{\includegraphics{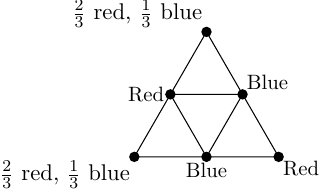}}
\label{f:hajos}
\caption{An optimal $2$-coloring of the Haj\'os graph}
\end{figure}

Consider, for example, the Haj\'os graph.  Figure~\ref{f:hajos} gives a $2$-coloring
of this graph
with defect~$4/3$ (and this is best possible in that any $2$-coloring has at least this much defect).
For another example, consider the complete graph on $3$ vertices.
Any $2$-coloring of~$K_3$
has defect at least $1$, but there are multiple optimal colorings:  color one
vertex red, one vertex blue, and the third vertex any combination of red
and blue.

Our objective is to minimize the defect of the graph.
Specifically, for a given number of
colors, what is the minimum defect that can be
obtained? If the number of colors is the chromatic
number, then of course the defect can be zero. But
if the number of colors is smaller, then there is a positive defect.

In the rest of the paper we proceed as follows.
In Section~\ref{s:prelim} we introduce some notation and
prove a few general results.
Thereafter in Section~\ref{s:twoColor}
we consider calculating the optimal $2$-colorings for several
graph families, including fans, wheels, complete multipartite graphs,
rooks graphs, and regular graphs. We give exact values in some cases and
bounds in others. We also pose several conjectures.
Finally in Section~\ref{s:complex} we observe that the decision problem is NP-hard.

\section{Preliminaries} \label{s:prelim}

Consider a graph $G$ and $k$ a number of colors.
For color $j$, let $f_j(v)$ be the usage of color $j$
on vertex $v$. Then the \dd{defect} of vertex $v$ is given by
\begin{equation}  \label{eq:defectV}
  \DEF (v) =  \sum_{w \in N(v)}   \sum_{j=1}^k  \min\left(   f_j(v), f_j(w) \right) .
\end{equation}
In general, the problem is to minimize
\[\
\max\limits_v  \DEF (v)
\]
over all colorings such that $f_j(v)$ is nonnegative
and $\sum_{j=1}^k f_j(v)=1$ for all vertices~$v$.  We denote this
minimum by $D(G,k)$, and call it the \dd{minimum defect}.

Note that
the existence of the minimum is guaranteed, since
the objective function is continuous and
the feasible region is a closed bounded set.
Further, the calculation is at least finite, since, for
example, we can prescribe which of $f_j(v)$ and
$f_j(w)$ are smaller in every min in Equation~\ref{eq:defectV} for each vertex $v$,
and thus $D(G,k)$ is the
minimum over exponentially many linear programs.

A related parameter is the \dd{minimum total defect}: we define
$\TD(G,k)$ to be the minimum over all colorings of
$\sum_{v\in V}  \DEF (v)$. However, here fractional colorings
do not yield anything new:

\begin{LEM}
For graph $G$ and number of colors $k$, there is a
$k$-coloring that achieves $\TD(G,k)$  in which every vertex
is monochromatic.
\end{LEM}
\begin{Proof}
Consider any vertex $v$ that is not monochromatic: say $f_1(v), f_2(v)>0$ with
$f_1(v)+f_2(v)=A$.
Then consider adjusting the coloring such that $f_1(v)=x$ and $f_2(v)=A-x$.
As a function of $x$, the defect of $v$ with any neighbor $w$ is a (piecewise-linear)
concave-down function. Thus, the total defect of the graph, as a function of~$x$,
is a concave-down function, and so its minimum is attained at an endpoint.
This means that one can either replace color~$1$ by color~$2$ or replace color~$2$ by
color~$1$ at $v$ without increasing the total defect. Repeated application of this replacement
yields a coloring with every vertex monochromatic.
\end{Proof}

For example, we will often use the fact that:

\begin{COR}  \label{c:totalComplete}
$\TD ( K_n, k ) = \lfloor n/k \rfloor (2n-k-\lfloor n/k \rfloor k)$.
\end{COR}

There are several fundamental results about monochromatic vertices for
minimum defect. One is that we may assume that there is
a monochromatic vertex of each color.

\begin{LEM} \label{l:monochromaticNodes}
Let $k$ be an integer and $G$ be a graph with at least $k$ vertices.
Then there is a $k$-coloring that achieves $D(G,k)$  that
has at least one monochromatic vertex for each color.
\end{LEM}
\begin{Proof}
Consider any optimal $k$-coloring, and consider each
color $j=1, 2, \ldots, k$ in turn. Each time, define
vertex $v_j$ as any vertex other than $v_1, \ldots, v_{j-1}$
with the largest usage of color~$j$; then
 recolor $v_j$ (if needed)
such that $f_j(v_j)=1$ and $f_i(v_j)=0$ for all $i\neq j$.
Such a recoloring does not increase the defect at any vertex.
So we will reach an optimal coloring with the desired property.
\end{Proof}

We next show that the minimum defect is
either $0$ or at least $1$.

\begin{LEM} \label{l:mindefect1}
For any graph $G$ and positive integer $k$, if
$D(G,k)>0$ then $D(G,k)\ge 1$.
\end{LEM}
\begin{Proof}
If every vertex is monochromatic, then the defect
is an integer and so the result is immediate. So consider any vertex $v$
that is not monochromatic.
If for any color $j$ we have $f_j(v) \ge f_j(w)$
for all neighbors $w$ of $v$,
then we can recolor $v$ to be monochromatically color $j$ without
increasing the defect of any vertex. So we may assume
that  for every
color $j$ at $v$, vertex $v$ has a neighbor~$w_j$ with $f_j(w_j) \ge f_j(v)$.
It follows that
\[
   \DEF (v) =  \sum_{w \in N(w)}   \sum_{j=1}^k \min\left(   f_j(v), f_j(w) \right)
   \ge  \sum_{j=1}^k \min\left(   f_j(v), f_j(w_j) \right)
  =  \sum_{j=1}^k   f_j(v) = 1 .
\]
The result follows.
\end{Proof}

For example, it follows from Lemma~\ref{l:mindefect1}
that the optimal $2$-coloring of the odd cycle $C_n$
has defect $1$.

\newpage

\begin{PROP} \label{p:complete}
The complete graph
$K_n$ has $D(K_n,k) =
\lceil n/k \rceil - 1$.
\end{PROP}
\begin{Proof}
This defect is achieved by (inter alia) coloring
each vertex with a single color and using each
color as equitably as possible. (This is trivially
the best coloring for total defect.)

To see that $\lceil n/k \rceil - 1$ is best possible, we proceed by induction on $n$,
noting that the result is trivial if $n \le k$. So assume
$n>k$. By Lemma~\ref{l:monochromaticNodes},
there is an optimal $k$-coloring with for each $j$
a vertex $v_j$ that is monochromatically color $j$.
Let $A= \{ v_1, \ldots, v_k\}$.
Then the defect of any other vertex $w$ in $G$ equals $1$ plus
the defect of $w$ in $G-A$. By the induction hypothesis, there
 exists a vertex in $G-A$ that has defect at least
$\lceil (n-k)/k \rceil -1 $ in $G-A$. This proves the lower bound.
\end{Proof}

\section{Two-Colorings of Some Graph Families}
\label{s:twoColor}

We now consider results for $2$-colorings. Unless otherwise
specified, we assume the colors are red and blue,
and denote the red usage at vertex $v$ by $r(v)$ (so that
the blue usage is $1-r(v)$). Unfortunately we are only able to provide
results for very specific families of graphs.

\subsection{Fans}

The \dd{fan}, denoted by $F_n$, is the graph obtained
from a path of order~$n$ by adding a new vertex and joining
it to every vertex of the path.

\begin{LEM} \label{l:F3}
In any $2$-coloring of $F_3$ it holds that $ \DEF (v)+ \DEF (w)\ge 2$ where $v$ and $w$
are the dominating vertices.
\end{LEM}
\begin{Proof}
Suppose the dominating vertices are $v$ and $w$ and the other two vertices
are $a$ and $b$.
Let $e_{xy}$ denote the overlap $\min( r(x), r(y) ) + \min( 1-r(x), 1-r(y) )$
between vertices $x$ and $y$.
Then $ \DEF (v)+ \DEF (w) = e_{va}+e_{vb}+e_{wa}+e_{wb}+2e_{vw}$;
further, because a triangle has total defect at least $2$ (Corollary~\ref{c:totalComplete}),
we have
$e_{va}+e_{wa}+e_{vw} \ge 1$ and $e_{vb}+e_{wb}+e_{vw} \ge 1$. The result follows.
\end{Proof}

Note that $F_1$ is just $K_2$ and $F_2$ is just $K_3$,
and so it holds that $D(F_1, 2)=0$ and $D(F_2, 2)=1$.
For the general cases of $F_n$, we have the following:

\begin{PROP} \label{p:Fan}
The minimum defect in a $2$-coloring of $F_n$ ($n\ge 3$) is
\[
D(F_n,2) = \frac{2 \lfloor n/3 \rfloor}{\lfloor n/3 \rfloor + 1} .
\]
\end{PROP}
\begin{Proof}
Say the path is $v_1v_2\ldots v_n$ with dominating vertex $h$.

We prove the upper bound by the following construction.
Set $x = 2/(\lfloor n/3 \rfloor + 1)$.
Let $r(h)=1$.
Let $r(v_i) = x$ if $i$ is a multiple of $3$, and $0$ otherwise.
It can readily be checked that every vertex $v_i$ has defect at most $2-x$, and
that vertex $h$ has defect $\lfloor n/3 \rfloor x$.
The result follows since both these values equal the claimed upper bound.

To prove the lower bound, it suffices to show that $D(F_n,2)\ge 2n/(n+3)$
if $n$ is a multiple of $3$. We partition the path $P_n$ into $n/3$ copies
of $P_3$; thus each $P_3$ along with vertex $h$ forms a
copy of $F_3$. Note that for each $1\le i \le n/3$, vertices $h$ and $v_{3i-1}$ are
the dominating vertices of the $i$\textsuperscript{th} copy of $F_3$. It follows from
Lemma~\ref{l:F3} that each copy of $F_3$ contributes at least $2$ to the sum
of the defects of $h$ and $v_{3i-1}$.
Therefore, $ \DEF (h) + \sum_{i=1}^{n/3}  \DEF (v_{3i-1}) \ge 2n/3$, whence
the result.
\end{Proof}

Note that the defect $D(F_n,2)$ tends to $2$ as $n$ increases. The fan is outerplanar.
Several researchers~\cite{BM85generalized,Mihok83partition} showed that one can
ordinarily $2$-color an outerplanar graph with defect at most $2$. However, we
conjecture that this bound can be improved slightly in the following sense:

\begin{CONJ} \label{c:Outerplanar}
$D(G,2)<2$ for any outerplanar graph $G$.
\end{CONJ}

\subsection{Wheels}

The \dd{wheel}, denoted by $W_n$, is the graph formed
from a cycle of order $n$ by adding a new vertex and
joining it to every vertex of the cycle.
The vertex of degree $n$ is called the \dd{hub} of the wheel.

\begin{PROP} \label{p:Wheel}
For $n\ge 3$, the minimum defect in a $2$-coloring of $W_n$ is
\[
D(W_n,2) = \frac{2 \lceil n/3 \rceil}{\lceil n/3 \rceil + 1} .
\]
\end{PROP}
\begin{Proof}
Let $x = 2/( \lceil n/3 \rceil + 1 )$ and let $D$ be a minimum independent
dominating set of the cycle. For a vertex $v$ on the cycle, let
$r(v)=x$ if $v\in D$, and $r(v)=0$ otherwise.
Let $r(h)=1$ for the hub~$h$.
It can readily be checked that every vertex on the cycle has defect at most $2-x$,
and that the hub has defect $x|D|$. The upper bound follows, since
$2-x = x|D| = 2 \lceil n/3 \rceil/( \lceil n/3 \rceil + 1 )$.

Next we prove the lower bound.
When $n=3k$, the lower bound follows directly from Proposition~\ref{p:Fan}.
Indeed, $D(W_{3k},2)\ge D(F_{3k},2) = 2k/(k+1)$. So we need to
establish the lower bound for $n=3k+1$ and $n=3k+2$.

Consider an optimal coloring of $W_{n}$ with hub $h$ and
cycle $v_1, v_2, \ldots, v_{n}, v_1$.
By Lemma~\ref{l:monochromaticNodes}, we may assume there exist vertices
$u$ and $u'$ such that $r(u)=0$ and $r(u')=1$.  There are two cases.

(a) \emph{Assume $h\notin \{u,u'\}$.}
Then we can form $k-1$ edge-disjoint copies of $P_3$ without using
vertex $h$, $u$, or~$u'$.
Let $S$ denote the set of centers of these copies.
By Lemma~\ref{l:F3}, it follows that the total defect of $S\cup\{ h \}$
within these copies is at least $2(k-1)$. Further, vertices $u$ and $u'$
together contribute defect $1$ to the hub $h$.
It follows that, in the graph as a whole, $ \DEF (h) +\sum_{s\in S}  \DEF (s) \ge 2k-1$, and so
$G$ has defect at least $(2k-1)/k$.
If $k\ge 2$, then $(2k-1)/k \ge (2k+2)/(k+2)$, and we are done.

So consider the case when $k=1$. Assume first that $n=5$.
Suppose $u$ and $u'$ are consecutive on the cycle; say $u=v_1$ and $u'=v_2$.
Then $G-\{u,u'\}$ is a copy of $F_3$. Since $u$ and $u'$ together contribute defect $1$
to $h$, it follows from Lemma~\ref{l:F3} that $ \DEF (h)+ \DEF (v_4) \ge 3$, and so $G$ has
defect at least $3/2$.
So assume without loss of generality that $u=v_1$ and $u'=v_3$.
If the hub has both two neighbors at least as red and
two neighbors at most as red, then it has defect at least $2$.
So without loss of generality, we may assume that $r(v_2), r(v_4), r(v_5) \le r(h)$.
Then, the defect that $h$ receives from $\{v_2, v_5\}$ is $2-2r(h)+r(v_2)+r(v_5)$,
and the defect that $u$ receives from $\{v_2, v_5\}$ is $2-r(v_2)-r(v_5)$.
That is, the sum of the defects that $h$ and $u$ receive from $\{v_2,v_5\}$ is at least $2$.
Since $h$ also receives defect~$1$ from $u$ and $u'$, it follows that $ \DEF (u)+ \DEF (h)\ge 3$,
and the result follows. The argument for $n=4$ is similar and omitted.

(b) \emph{Assume $h\in \{u,u'\}$.}
Say $u=v_1$ and $u'=h$.  Consider $v_n$. It receives
defect $1$ from $\{v_1,h\}$. Let index $j$ be such that $v_j$ is the
redder vertex of $v_{n-1}$ and $v_{n}$. Then $v_{n-1}$ and $v_n$ have at least
$1-r(v_j)$ of blue overlap and so $ \DEF (v_n) \ge 2 - r(v_j)$.

Further, one can form $k$ edge-disjoint
copies of $P_3$ without using vertex $h$ or $v_j$. Let $S$ denote the set of centers of these copies.
By Lemma~\ref{l:F3} and noting that the hub $h$ receives defect $r(v_j)$ from vertex $v_j$,
it follows that $ \DEF (h)+\sum_{s\in S}  \DEF (s) \ge 2k + r(v_j)$.

Thus $ \DEF (h)+ \DEF (v_n) + \sum_{s\in S}  \DEF (s) \ge 2k+2$, whence the result.
\end{Proof}

\subsection{Complete multipartite graphs and compositions}

We consider here complete multipartite graphs. These can
be thought of as taking a complete graph and replacing
each vertex by an independent set with the same adjacency.
In general, we define $G[aK_1]$ to be the \dd{composition} of
$G$ with the empty graph on $a$ vertices; that is, the
graph obtained by replacing every vertex $v$ of $G$
with a set~$I_v$ of size $a$ such that a vertex of
$I_v$ is adjacent to a vertex of $I_w$ if and only if
$v$ and $w$ are adjacent in $G$.

There are two simple bounds:

\begin{PROP} \label{p:expansion} For any graph $G$, \\
(a) $\TD(G[aK_1], k ) \ge a^2 \TD(G,k)$.  \\
(b) $D(G[aK_1],k) \le aD(G,k)$.
\end{PROP}
\begin{Proof}
(a) The bound follows by applying the lower bound to the $a^n$ copies of $G$
and averaging. \\
(b) Take the optimal coloring of $G$ and replicate it: give every vertex of $I_v$ the
color of vertex $v$.
\end{Proof}

We let $\MULTI{m}{a}$ denote the
complete $m$-partite graph with $a$ vertices in each partite set;
that is $\MULTI{m}{a}=K_m[aK_1]$.
It follows that:

\begin{PROP} \label{p:CMspecialMultiple}
If $m$ is a multiple of $k$, then the complete multipartite graph
$\MULTI{m}{a}$ can be $k$-colored with defect $(m/k -1)a$,
and this is best possible.
\end{PROP}

But if $m$ is not a multiple of $k$, the result is not
clear. Perhaps the following is true.

\begin{CONJ} \label{conj:compMulti}
The minimum defect in a $k$-coloring of
$\MULTI{m}{a}$ is $(\lceil m/k \rceil -1)a$.
\end{CONJ}

In fact, we do not have an example that precludes it being the
case that it always holds that $D(G[aK_1],k) = aD(G,k)$.

We shall prove Conjecture~\ref{conj:compMulti} for $2$
colors.
We need the following definitions. Define a vertex $x$ as
\dd{large} if $r(x) > 1/2$ and
\dd{small} if $r(x) < 1/2$.
Also we let $N(x)$ denote the set of neighbors of $x$,
$U(x)$ denote the set of vertices $y$ in $N(x)$ with $r(y) \ge r(x)$, and
$L(x)$ denote the set of vertices $y$ in $N(x)$ with $r(y) < r(x)$.

We also need the following observations and lemmas. Some of them are
very easy to verify and so the proofs are omitted:

\begin{OBS} \label{p:MiddleVertex}
If $r(x) = 1/2$, then $\DEF(x) \ge |N(x)|/2$.
\end{OBS}

\begin{OBS} \label{p:SameParity}
If two vertices are both large (or both small), then the overlap
between them is greater than $1/2$.
\end{OBS}

\begin{OBS} \label{p:UL}
$\DEF(x) \ge \min\left( |U(x)|, |L(x)| \right)$.
\end{OBS}

\begin{LEM} \label{l:LargeSmallVertex}
$\DEF(x) \ge |N(x)|/2$ if either \\
(a) $x$ is large and $|U(x)|\ge |L(x)|$, \\
or (b) $x$ is small and $|U(x)|\le |L(x)|$.
\end{LEM}
\begin{Proof}
It suffices to prove it for the case that $x$ is large. We pair
each vertex in $L(x)$ with a vertex in $U(x)$.
Then each pair contributes at least $1$ to $\DEF(x)$.
By Observation~\ref{p:SameParity}, each of the remaining vertices in
$U(x)$ contributes more than $1/2$ to $\DEF(x)$.
Hence $\DEF(x) \ge |N(x)|/2$.
\end{Proof}

\begin{LEM} \label{l:DifferentParity}
If $x$ is large and $y$ is small, then
$\max \left( \DEF(x), \DEF(y) \right) \ge |N(x) \cap N(y)|/2$.
\end{LEM}
\begin{Proof}
If $|U(x)|\ge |L(x)|$, then we have $\DEF(x) \ge |N(x)|/2 \ge |N(x) \cap N(y)|/2$
by Lemma~\ref{l:LargeSmallVertex}. So we may assume
$|U(x)| < |L(x)|$. Similarly we may assume $|U(y)| > |L(y)|$.
Note that we can increase $r(x)$ to $1$ and decrease $r(y)$ to $0$ without increasing
the defect of either vertex. It follows that $\DEF(x) + \DEF(y)$
is at least their common degree, whence the result.
\end{Proof}

\begin{LEM} \label{l:ShiftToExtrema}
If all neighbors of $x$ are large (small), then $r(x)$ can be changed to $0$ $(1)$
without increasing the defect of any vertex.
\end{LEM}
\begin{Proof}
It suffices to prove it for the case that all neighbors of $x$ are large. Let $v$ be any neighbor of $x$.
The overlap between them is $1-|r(v)-r(x)|$. If $r(x)$ is changed to $0$, then
the overlap becomes $1-r(v)$. Since $r(v) > 1/2$, we have
$1-|r(v)-r(x)| \ge 1-r(v)$ and the conclusion follows.
\end{Proof}

\begin{PROP} \label{p:CMspecial}
The minimum defect in a $2$-coloring of $\MULTI{m}{a}$ is $(\lceil m/2 \rceil-1)a$.
\end{PROP}
\begin{Proof}
Such defect
is attained by coloring all vertices in $\lfloor m/2 \rfloor$
of the partite sets with red, and the remaining vertices blue.
So we need to prove that this is best possible.

If $m$ is even, Proposition~\ref{p:expansion}
shows that $\TD(\MULTI{m}{a}, 2) \ge m( m/2-1)a^2$,
and thus some vertex has defect at least
$(m/2-1)a$. So assume $m$ is odd.

If there is a vertex $v$ in the graph with $r(v) = 1/2$, then
$\DEF (v) \ge (m-1)a/2 = (\lceil m/2 \rceil-1)a$ by Observation~\ref{p:MiddleVertex}.
Also, if there is a partite set that contains both a large
vertex and a small vertex, then we are okay by Observation~\ref{l:DifferentParity}.

Hence we may assume every partite set contains either only
large vertices or only small vertices. Without loss of generality,
assume at least $(m+1)/2$ partite sets contain only large vertices.
Let $x$ be the large vertex with \emph{minimum} $r(x)$.
Note that $|U(x)| \ge (m-1)a/2 \ge |L(x)|$, and therefore
$\DEF (x) \ge (m-1)a/2 = (\lceil m/2 \rceil-1)a$ by Observation~\ref{l:LargeSmallVertex}.
\end{Proof}

\begin{PROP} \label{p:CMabc}
The minimum defect in a $2$-coloring of the complete tripartite graph $K_{a,b,c}$
with $a \le b \le c$ is $bc/(b+c-a)$.
\end{PROP}
\begin{Proof}
Let $A$, $B$, and $C$ denote the partite sets of order $a$, $b$, and $c$, respectively.
The upper bound is attained by coloring
all vertices $v$ in $A$ with $r(v)=0$, all vertices in $C$ with $r(v)=1$, and
all vertices in $B$ with $r(v)=x$, where $x$ is chosen to give the vertices in $A$ and
$B$ the same defect, namely $x = (b-a)/(b+c-a)$.

Now we prove the lower bound.
Let $x_1, x_2, \ldots, x_a$ be the vertices in $A$ with
$r(x_1) \le r(x_2) \le \ldots \le r(x_a)$,
$y_1, y_2, \ldots, y_b$ be the vertices in $B$ with
$r(y_1) \le r(y_2) \le \ldots \le r(y_b)$,
and $z_1, z_2, \ldots, z_c$ be the vertices in $C$ with
$r(z_1) \le r(z_2) \le \ldots \le r(z_c)$.
There are two possible cases.

\emph{\textbf{Case 1}: $a \le b \le c \le a+b$.}

Then we have $(b+c)/2 \ge (a+c)/2 \ge (a+b)/2 \ge bc/(b+c-a)$.
If there is a vertex $v$ in the graph with $r(v) = 1/2$, then
the conclusion follows from Observation~\ref{p:MiddleVertex}. Also, if
there is a partite set that contains both a large
vertex and a small vertex, then the conclusion follows from
Lemma~\ref{l:DifferentParity}.
Hence we may assume that every partite set contains either only
large vertices or only small vertices, and by symmetry we only need to
consider the following four subcases:

\emph{\textbf{Case 1.1}: all vertices in the graph are large.}

By Observation~\ref{p:SameParity},
we have $\DEF(x_i) \ge (b+c)/2$ for every $1 \le i \le a$. So the conclusion follows.

\emph{\textbf{Case 1.2}: all vertices in $A$ are small and all the other vertices are large.}

Let $u$ be the large vertex with \emph{minimum} $r(u)$.
By Lemma~\ref{l:LargeSmallVertex}, $\DEF(u) \ge (a+b)/2$ and the conclusion follows.

\emph{\textbf{Case 1.3}: all vertices in $B$ are small and all the other vertices are large.}

By Lemma~\ref{l:ShiftToExtrema}, we may assume $r(y_j)=0$ for every $1 \le j \le b$.
If $r(x_1) \le r(z_1)$, then by Lemma~\ref{l:LargeSmallVertex},
$\DEF(x_1) \ge (b+c)/2$.
So assume $r(x_1) > r(z_1)$. We have
\begin{eqnarray*}
\DEF(x_a)   &=& b(1-r(x_a)) + \sum\limits_{k=1}^{c}(1-|r(x_a)-r(z_k)|) \\
        &\ge& b(1-r(x_a)) + \sum\limits_{k=1}^{c}(r(x_a)+r(z_k)-1) \\
        &=& (b-c)(1-r(x_a)) + \sum\limits_{k=1}^{c}r(z_k) \\
        &\ge& (b-c)(1-r(x_a)) + c~r(z_1),
\end{eqnarray*}
and
\begin{eqnarray*}
\DEF(z_1) &=& \sum\limits_{i=1}^{a}(1-(r(x_i)-r(z_1))) + b(1-r(z_1)) \\
       &=& \sum\limits_{i=1}^{a}(1-r(x_i)) + (a-b)r(z_1) + b \\
       &\ge& a(1-r(x_a)) + (a-b)r(z_1) + b.
\end{eqnarray*}

Hence, $(b-a)\DEF(x_a) + c\DEF(z_1) \ge [(b-a)(b-c)+ac](1-r(x_a)) + bc \ge bc$.
It follows that $\max \left( \DEF(x_a), \DEF(z_1) \right) \ge bc/(b+c-a)$.

\emph{\textbf{Case 1.4}: all vertices in $C$ are small and all the other vertices are large.}

By Lemma~\ref{l:ShiftToExtrema}, we may assume $r(z_k)=0$ for every $1 \le k \le c$.
If $r(x_1) \le r(y_1)$, then we have
\begin{eqnarray*}
\DEF(x_1) &=& c(1-r(x_1)) + \sum\limits_{j=1}^{b}(1-(r(y_j)-r(x_1))) \\
       &=& c+(b-c)r(x_1)+\sum\limits_{j=1}^{b}(1-r(y_j)) \\
       &\ge& c+(b-c)r(x_1),
\end{eqnarray*}
and
\begin{eqnarray*}
\DEF(y_1) &=& \sum\limits_{i=1}^{a}(1-|r(x_i)-r(y_1)|) + c(1-r(y_1)) \\
       &\ge& \sum\limits_{i=1}^{a}(r(x_i)+r(y_1)-1) + c(1-r(y_1)) \\
       &=& (c-a)(1-r(y_1)) + \sum\limits_{i=1}^{a}r(x_i) \\
       &\ge& \sum\limits_{i=1}^{a}r(x_i) \\
       &\ge& a~r(x_1).
\end{eqnarray*}

Hence, we have
\begin{eqnarray*}
b\DEF(x_1)+(c-a)\DEF(y_1) &\ge& bc+[b(b-c)+(c-a)a]r(x_1) \\
&=& bc+(b+a-c)(b-a)r(x_1) \\
&\ge& bc.
\end{eqnarray*}

It follows that $\max \left( \DEF(x_1), \DEF(y_1) \right) \ge bc/(b+c-a)$.

Similarly, if $r(x_1) > r(y_1)$, then it can be verified that
\begin{eqnarray*}
\DEF(x_1) \ge (c-b)(1-x_1) + \sum\limits_{j=1}^{b}r(y_j) \ge b~r(y_1),
\end{eqnarray*}
and
\begin{eqnarray*}
\DEF(y_1) = \sum\limits_{i=1}^{a}(1-r(x_i)) + c + (a-c)r(y_1) \ge c + (a-c)r(y_1).
\end{eqnarray*}

Hence, we have
\[
(c-a)\DEF(x_1)+b\DEF(y_1) \ge b(c-a)r(y_1) + bc + b(a-c)r(y_1) = bc.
\]

It follows that $\max \left( \DEF(x_1), \DEF(y_1) \right) \ge bc/(b+c-a)$.

\emph{\textbf{Case 2}: $a \le b < a+b < c$.}

Then we have $(b+c)/2 \ge (a+c)/2 > c/2 > bc/(b+c-a)$.
By Observation~\ref{p:MiddleVertex} and Lemma~\ref{l:DifferentParity},
we only consider the case that the vertices of $A \cup B$ are
either all large or all small.
Without loss of generality, assume that they are all large.
Then by Lemma~\ref{l:ShiftToExtrema}, we may assume that
$r(z_k)=0$ for every $1 \le k \le c$.

If $r(x_1) \le r(y_1)$, then $\DEF(x_1) \ge b$ by Observation~\ref{p:UL}.
So assume $r(x_1) > r(y_1)$. But then by the same argument as that in Case 1.4,
we have $\max \left( \DEF(x_1), \DEF(y_1) \right) \ge bc/(b+c-a)$.
\end{Proof}

For another composition, consider $C_m[aK_1]$ where $m$ is odd.
We now prove that $D( C_m[2K_1], 2 ) = 2$. There are at least
two different optimal colorings.
The first such coloring is obtained by taking
an optimal coloring for~$C_m$ and replicating it. The second such
coloring is obtained by, for each copy of $2K_1$, coloring
one vertex red and one vertex blue.

\begin{PROP} \label{p:cycleSpecial}
For $m$ odd, $D(C_m[2K_1], 2) = 2$.
\end{PROP}
\begin{Proof}
Consider a $2$-coloring of $C_m[2K_1]$.
We need to show that the defect is at least~$2$.
As in the proof of Proposition~\ref{p:CMspecial}, we may
assume that every copy of $2K_1$ contains either two large
vertices or two small vertices. Since $m$ is odd, it follows that
there must be two adjacent copies of the same type.
Without loss of generality, assume $u_1$ and $u_2$ are adjacent to
$v_1$ and $v_2$ with all four vertices being large.
If any $x\in \{u_1,u_2,v_1,v_2\}$ has $|U(x)|\ge 2$,
then the lower bound follows from Lemma~\ref{l:LargeSmallVertex}(a).
Therefore we may assume that $|U(x)|\le 1$ for
every $x\in \{u_1,u_2,v_1,v_2\}$. This means that each $u_i$ is
redder than some $v_j$ and vice versa, a contradiction.
\end{Proof}

\subsection{Rooks graphs and Cartesian products}

Recall that the \dd{Cartesian product} of graphs $G$ and $H$, denoted $G \Box H$,
is the graph whose vertex set is $V(G) \times V(H)$, in which two vertices
$(u_1, u_2)$ and $(v_1, v_2)$ are adjacent if $u_1v_1 \in E(G)$ and $u_2=v_2$,
or $u_1=v_1$ and $u_2v_2 \in E(H)$.

We will need the obvious lower bound for the total defect of Cartesian products.

\begin{PROP} \label{p:cart}
Let $G$ and $H$ be graphs of order $m$ and $n$ respectively.
Then \\ $\TD(G \Box H, k) \ge m \TD(H,k) + n \TD(G,k)$.
\end{PROP}
\begin{Proof}
The defect of a vertex in the product is the sum of the defects in its
copies of $G$ and $H$.
\end{Proof}

The \dd{rooks graph}, denoted by $K_m \Box K_n$, is the
Cartesian product of the complete graphs $K_m$ and $K_n$.  We will
denote the vertices by
 $(i, j)$ with $1 \le i \le m, 1 \le j \le n$.

\begin{LEM} \label{l:rooksUB}
The rooks graph $K_m \Box K_n$ can be $2$-colored with defect
$\lceil m/2 \rceil + \lceil n/2 \rceil - 2$.
\end{LEM}
\begin{Proof}
Color vertex $(i, j)$ with red if $i$ and $j$ have the same parity
and blue otherwise.
\end{Proof}

\begin{COR} \label{p:rooksEE}
Let $m$ and $n$ be even integers. Then
$D(K_m \Box K_n, 2) = m/2 + n/2  - 2$.
\end{COR}
\begin{Proof}
The upper bound follows from Lemma~\ref{l:rooksUB}.
The lower bound follows from Proposition~\ref{p:cart}, since
$\TD( K_s , 2) = s(s/2-1)$ for $s$ even (Corollary~\ref{c:totalComplete}), and thus
$\TD(K_m \Box K_n, 2) \ge mn(n/2-1)+nm(m/2-1)$.
\end{Proof}

We show below that the upper bound in Lemma~\ref{l:rooksUB}
is not always optimal. In fact we conjecture that it is never
optimal when $m$ and $n$ are both odd, except for the case that
$m=n=3$.

\begin{LEM}
$D(K_3 \Box K_3, 2 )=2$.
\end{LEM}
\begin{Proof}
The upper bound is from Lemma~\ref{l:rooksUB}.

We have two proofs of the lower bound, one by computer and one
by hand. Both proofs entail converting the question to
a set of linear programs.

Observe that given a coloring, one can generate an acyclic orientation
by orienting each edge from smaller to larger proportion of red (with ties broken
by vertex number say). Further, if $N_1$ is the set of neighbors of vertex $v$ with more red
and $N_2$ is the set of neighbors of $v$ with less red, then
Equation~\ref{eq:defectV} simplifies to
\[
 \DEF (v) =  |N_1| r(v) + |N_2| b(v) +
\sum_{w \in N_2}  r(w) + \sum_{w \in N_1}  b(w) ,
\]
where $b(x)=1-r(x)$.

So, the proof is to enumerate the acyclic orientations.
For each such orientation, we add the constraints that
$r(u) \le r(v)$ for all arcs $uv$. That is, minimizing
the defect for a given orientation is a linear program.

Further, if any
vertex has in- and out-degree~$2$ for the orientation, the defect is
definitely at least~$2$ (by Observation~\ref{p:UL}).
With several
pages of calculation or by using a computer, one can show that $K_3 \Box K_3$ has eight
acyclic orientations (up to symmetry) that need to be considered, and then solve the eight associated
linear programs. We omit the details.
\end{Proof}

In contrast, we found a
coloring of $K_3 \Box K_5$ that beats the bound of Lemma~\ref{l:rooksUB}:

\begin{LEM} \label{l:rooks35}
$D(K_3 \Box K_5, 2 ) \le  38/13$.
\end{LEM}
\begin{Proof}
A $2$-coloring of $K_3 \Box K_5$ is shown below. The element $(i, j)$
of the matrix is the red-usage on vertex $(i, j)$.
\[
\left[ \begin{array}{ccc}
          0      & 8/13    & 0 \\
          0      & 0       & 8/13 \\
          1      & 11/13   & 0 \\
          1      & 0       & 11/13 \\
          6/13   & 1       & 1
\end{array} \right]
\]
It can be verified that the defect of the coloring
is $38/13$.
\end{Proof}

The above coloring can be extended to show that
Lemma~\ref{l:rooksUB} is not optimal for
$m=3$ and $n$ odd, $n \ge 5$,
and indeed that $D( K_3 \Box K_n,2 )  \le
n/2+11/26$ in this case.
However, this is still not best possible.
For example, one can get defect
$42/11$ for $K_3 \Box K_7$
and defect $14/3$ for $K_3 \Box K_9$ by the colorings illustrated:

\[
\left[
\begin {array}{ccc}
1&0&1\\
4/11&1&1 \\
0&8/11&0 \\
1&1&4/11 \\
1/11&0&1\\
0&8/11&0 \\
1&0&1/11
\end {array}
\right]
\qquad\qquad\qquad
 \left[
 \begin {array}{ccc}
 2/3&0&0\\
 0&2/3&0\\
 0&0&2/3\\
 0&1&1\\
 0&1&1\\
 1&0&1\\
 1&0&1\\
 1&1&0\\
 1&1&0
 \end {array} \right]
\]

We used simulated annealing computer search to find upper bounds. Though we
have no exact values, it seems to us that the heuristic computer search results suggest the following:

\begin{CONJ} \label{c:rooksConj1}
(a) If $m+n$ is odd, then $D(K_m \Box K_n,2) = (m+n-3)/2$. \\
(b) If $m \ge 3$, $n \ge 3$, and $mn$ is odd and greater than $9$, then
$D(K_m \Box K_n,2) < \lceil m/2 \rceil + \lceil n/2 \rceil - 2$.
\end{CONJ}

Note that Conjecture~\ref{c:rooksConj1} $(a)$ is trivially true for the case
that $m=2$ (or $n=2$), since $D(K_2 \Box K_n,2) \ge D(K_n,2)= (n-1)/2$.

Proposition~\ref{p:cart} yields the following lower bounds:

\begin{COR} \label{c:rooksOldLB} \ \\
(a) If both $m$ and $n$ are odd, $D(K_m \Box K_n, 2) \ge (m+n)/2-2+1/(2m)+1/(2n)$. \\
(b) If $m$ is even and $n$ is odd, $D(K_m \Box K_n, 2) \ge (m+n)/2-2+1/(2n)$.
\end{COR}

For more colors we have one trivial observation:
that
$D(K_n \Box G, k) = \lceil n/k \rceil -1$ for any $k$-partite graph $G$, as a corollary of
Proposition~\ref{p:complete}.

\subsection{Regular graphs}

Lov\'asz~\cite{Lovasz66decomposition} showed that we can ordinarily $2$-color a cubic graph
with defect at most $1$. Therefore $D(G,2)=1$ for all
nonbipartite cubic graphs $G$.

For a $4$-regular graph, Lov\'asz's result shows
that one can ordinarily $2$-color it with defect at most $2$.
We conjecture that this can be improved.  Proposition~\ref{p:cycleSpecial}
shows that the composition $G=C_{m}[2K_1]$ where $m$ is odd has $D(G,2)=2$.
Using simulated annealing (that is, a randomized search for a coloring),
the computer can
find a $2$-coloring with defect smaller than $2$ for all $4$-regular graphs on up to
$14$ vertices, except for the compositions of odd cycles, and the two graphs
$K_5$  and $K_3\Box K_3$,  which we saw earlier have minimum defect $2$.
We conjecture a general behavior:

\begin{CONJ}
Apart from $G=C_{m}[2K_1]$ where $m$ is odd, it holds that $D(G)<2$ for
all but finitely many connected $4$-regular graphs.
\end{CONJ}

\section{Complexity} \label{s:complex}

Unsurprisingly, it is NP-hard to determine if there is a
coloring with defect at most some specified~$d$.

One way to see this is that fractional defect $2$-coloring is NP-hard even for $d=1$.
One can extend Lemma~\ref{l:monochromaticNodes} to show that in graphs of
minimum degree at least $3$, a $2$-coloring with defect $1$
can only be a coloring with monochromatic vertices.
Thus the fractional defect $2$-coloring problem is equivalent to the
ordinary defective $2$-coloring problem in such graphs. The latter problem was shown to
be NP-hard by Cowen~\cite{Cowen93complexity}. (Actually, we need ordinary $1$-defect coloring to
be NP-hard in graphs with minimum degree at least $3$.
But one can transform a graph to having minimum degree at
least $3$ without changing the coloring property by
adding, for each vertex $v$, a copy of $K_4$ and joining $v$ to one vertex of the $K_4$.)

\section{Acknowledgements} \label{s:Acknowledgements}

The authors would like to thank anonymous referees for helpful comments.
This work is part of the second author's Ph.D. Dissertation at
Clemson University~\cite{xu2016generalized}, and he wishes to thank
his doctoral committee members for their guidance and suggestions.

%\newpage

\bibliographystyle{plain}
\bibliography{fractional}

\end{document}